\documentclass[12pt]{amsart}

\input{xy}
\xyoption{all}

\usepackage{graphicx}
\usepackage{color}
\usepackage{amssymb}
\usepackage{amsmath}
\usepackage{graphics}
\usepackage{graphics}
\usepackage{enumerate}

\newtheorem{theorem}{Theorem} [section]
\newtheorem{prop}[theorem]{Proposition}

\newtheorem{cor}[theorem]{Corollary}

\numberwithin{equation}{section}
\numberwithin{figure}{section}

\newcommand\C{{\mathbb C}}

\renewcommand\P{{\mathbb P}}

\newcommand\N{{\mathbb N}}

\renewcommand\phi{\varphi}

\newcommand\Rat  {\mathrm{Rat}} %space of rational maps
\newcommand\Prep {\operatorname{PrePer}}
\newcommand\Per {\operatorname{Per}}

\newcommand\Crit {\operatorname{Crit}}

\begin{document}

\title{RATIONAL FUNCTIONS WITH IDENTICAL MEASURE OF MAXIMAL ENTROPY}

\author{Hexi Ye}

\begin{abstract}
We discuss when two rational functions $f$ and $g$ can have the same measure of maximal
entropy.  The polynomial case was completed by (Beardon, Levin, Baker-Eremenko, Schmidt-Steinmetz, etc.,
1980s-90s), and we address the rational case following Levin-Przytycki (1997).   We
show:   $\mu_f = \mu_g$ implies that  $f$ and $g$ share an iterate
($f^n = g^m$ for some $n$ and $m$) for general $f$ with degree $d
\geq 3$.  And for generic $f\in \Rat_{d\geq 3}$,  $\mu_f = \mu_g$ implies $g=f^n$ for some $n
\geq 1$.
For generic $f\in \Rat_2$, $\mu_f = \mu_g$ implies that $g= f^n$ or $\sigma_f\circ f^n$  for some $n\geq 1$, where $\sigma_f\in PSL_2(\C)$ permutes
two points in each fiber of $f$.
 Finally, we construct examples of $f$ and $g$ with $\mu_f = \mu_g$
such that $f^n \neq \sigma\circ g^m$ for any $\sigma \in PSL_2(\C)$ and $m,n\geq 1$.
\end{abstract}

\maketitle

% For editing purposes only
%\tableofcontents

% suppress page number 1
\thispagestyle{empty}

%%%%%%
%%%%%%

\section{Introduction}

Let $f:\P^1 \to \P^1$ be a rational function with degree $d_f\geq 2$, where the projective space $\P^1$ is defined over $\C$. There is a unique probability measure $\mu_f$ on $\P^1$, which is invariant under $f$ and has support equal to the Julia set $J_f$ of $f$, achieving maximal  entropy $h_{\mu_f}=\log d$ among all the $f$-invariant probability measures; see \cite{L3} and \cite{FLM}.

In this article, we study rational functions with the same measure of maximal entropy. It is well known that $\mu_f=\mu_{f^n}$ for all iterates $f^n$ of $f$, and commuting rational functions have common measure of maximal entropy. In the polynomial case, having the same measure of maximal entropy is equivalent to having the same Julia set. During the 1980s and 90s, pairs of polynomials with identical Julia set were characterized; see \cite{SS} \cite{B2} \cite{B3} and \cite{BE}. The strongest result is: given any Julia set $J$ of some non-exceptional polynomial,  there is a polynomial $p$, such that the set of all polynomials with Julia set $J$ is
\begin{equation}\label{mpoly}
\{\sigma \circ p^n|n\geq 1 \textup{ and } \sigma \in \textup{$\Sigma_J$}\},
\end{equation}
where $\Sigma_J$ is the set of complex affine maps on $\C$ preserving $J$. By definition, a rational function is {\bf exceptional} if it is conformally conjugate to either a power map, $\pm$Chebyshev polynomial, or a Latt$\grave{\textup{e}}$s map.
From (\ref{mpoly}), if $f$ and $g$ are two non-exceptional  polynomials with $\mu_f=\mu_g$, then there exists  $\sigma(z)=az+b$ preserving $\mu_f$ with
\begin{equation}\label{polycom}
f^n=\sigma \circ g^m \textup{ for some $m, n\geq 2$.}
\end{equation}

\medskip

However, unlike the polynomial case, there exist non-exceptional rational functions with the same maximal measure but not related by the formula (\ref{polycom}).
\begin{theorem}\label{counter}
There exist non-exceptional rational functions $f$ and $g$ with degrees $\geq 2$ and $\mu_f=\mu_g$, but
\begin{equation}\label{neqequ}
f^n\neq \sigma\circ g^m \textup{ for any $\sigma\in PSL_2(\C)$ and $n,m\geq1$}.
\end{equation}
Specifically, for $R, S, T$ being rational functions with degrees $\geq 2$ such that
\begin{itemize}
\item For any $\sigma\in PSL_2(\C)$, we have  $R\neq \sigma \circ S$.
\item $T\circ R=T\circ S$.
\end{itemize}
we set $f=R\circ T$ and $g=S\circ T$, then $\mu_f=\mu_g$ and they satisfy  (\ref{neqequ}).

%Let $f=\frac{(z^3-3z)^2+1}{z^3-3z}$ and $g=\frac{e^{2\pi i/3 }(z^3-3z)^2+e^{4\pi i/3}}{z^3-3z}$ with $a\neq 0$ and $\omega^2+\omega+1=0$. We have $\mu_f=\mu_g$ and none of them is exceptional.  For any $n,m\geq 1$ and any M$\ddot{\textup{o}}$bius transformation $\sigma$, we have
%$$f^n\neq \sigma\circ g^m.$$
\end{theorem}
\noindent The existence of the triples ($R,S,T$) in Theorem \ref{counter} is equivalent to the existence of an irreducible component of
$$V_T=\{(x,y):T(x)=T(y)\}\subset \P^1\times\P^1$$ with bidegree ($r,r$), $r\geq 2$, and normalization of genus $0$. Explicit examples for such triples  ($R,S,T$) are provided in Section 3.
\medskip

Let $\Rat_d$ be the set of all   rational functions with degree $d\geq 2$.  The space $\Rat_d$ sits inside $\P^{2d-1}(\C)$, and it is the complement of the zero locus of an irreducible homogenous polynomial (the resultant) on $\P^{2d-1}$; therefore  $\Rat_d$ is an affine variety.
For any rational function $f\in \Rat_d$, denote by $M_f$ the set of all  rational functions with the same maximal entropy measure as $f$. As we discussed before, when $f$ is non-exceptional and conjugate to some polynomial, $M_f$ has very simple expression as in (\ref{mpoly}) by Corollary \ref{poly id}.  However, from  Theorem \ref{counter}, we do not have the conclusion of (\ref{polycom}) for all non-exceptional rational functions $f$ and $g$ with $\mu_f=\mu_g$, even we replace $\sigma$ by any M$\ddot{\textup{o}}$bius transformation.
Levin \cite{L1} \cite{L2} showed $M_f\bigcap\Rat_n$ is a finite set unless $f$ is  conjugate to the power function $z^{\pm d}$.

%The known result is that it is a finite group unless $f(z)$ is exceptional, see \cite{??}.

\medskip
For convenience,
in the rest of this paper, {\bf generic} means with exception of at most countably many proper Zariski closed subsets; {\bf general} means with exception of  some proper Zariski closed subset.

In this article, we prove,
\begin{theorem}\label{generic rational}
Let $\Rat_d$ be the set of all rational functions with degree $d\geq 2$. For generic rational functions $f\in \Rat_d$, we have
\begin{itemize}
\item $M_f=\{f, f^2, f^3, \cdots\}$, when $d\geq 3$,
\item $M_f=\{f, \sigma_f\circ f, f^2, \sigma_f\circ f^2, f^3,\sigma_f\circ f^3, \cdots\}$, when $d=2$,
\end{itemize}
where $\sigma_f$ is the unique M$\ddot{\textup{o}}$bius transformation permuting the fibers of $f$.
\end{theorem}

The proof of Theorem \ref{generic rational} is mainly based on the next two theorems. The first, Theorem \ref{mainthm}, asserts that for general rational functions with degree $d\geq 3$, having the same measure of maximal entropy is the same as sharing an iterate. We will say a critical value of $f\in \Rat_d$ is {\bf simple} if its preimage contains exactly one critical point counted with multiplicity.

\begin{theorem}\label{mainthm}
Let $f$ be a rational function with degree $d_f\geq 3$, and $f$ has at least three simple  critical values. Then for any rational function $g$ with degree $d_g\geq 2$ and $\mu_f=\mu_g$, we have
$$f^n=g^m$$
for some integers $n,m\geq 1$.
%Let $f$ be a rational function with degree $d\geq 3$ and $2d-2$ distinct grand orbits for the critical points, moreover none of the critical points is preperiodic. Assume $g$ is a rational function with degree $\geq2$ and $\mu_g=\mu_f$. Then $f$ and $g$ have some common iterate, i.e.
%$$f^n=g^m, \textup{ for some $n, m\geq 2,$ and $n|m$.}$$
\end{theorem}

Theorem \ref{mainthm} only works for rational functions with degree $d\geq 3$. This is because for degree $d=2$ case, there is a special nontrivial symmetry $\sigma_f \in PSL_2(\C)$ for each
$f\in \Rat_2$, see \cite{LP}. The symmetry $\sigma_f$ is  the unique M$\ddot{o}$bius transformation permuting the fibers of $f$ with
$\sigma_f^2(z)=z$. As $\sigma_f$ permutes the points in the fiber of $f$, we have  $f\circ \sigma_f=f$ and then $\sigma_f$ preserves $\mu_f$.
 Hence we have $\mu_{\sigma_f\circ f}=\mu_f$. For any  $f\in \Rat_2$, let $g=\sigma_f\circ f$. It satisfies $\mu_g=\mu_f$ and $g^n=\sigma_f\circ
f^n\neq f^n$ for any $n\geq 1$. In other words, $f$ and $g$ have the same maximal measure but they never share an iterate. In all, for any $f\in \Rat_d$, we have obvious  relations: $\{f, f^2, f^3, \cdots\}\subset M_f$, and when $d=2$, $\{f, \sigma_f\circ f, f^2, \sigma_f\circ f^2, f^3,\sigma_f\circ f^3, \cdots\}\subset M_f$. Theorem \ref{generic rational} asserts that, generically, there is no other rational function in $M_f$. However, it is still not
known  whether we can replace ``generic" in Theorem \ref{generic rational} by ``general", which will greatly improve the result; at least it is clear from (\ref{mpoly}) that the statements in Theorem \ref{generic rational} are satisfied for general polynomials.

\medskip
Let  $d\geq 2$ and $ n\geq 1$ be integers. There is a  regular map between affine varieties:
$$\phi_{d,n}: \Rat_d \to \Rat_{d^n},$$
defined by $\phi_{d,n}(f)=f^n$. We call it the iteration map of rational functions.

The next result, Theorem \ref{injective}, states that the iteration map is one-to-one for general points.

\begin{theorem}\label{injective}
Let $\phi_{d,n}: \Rat_d \to \Rat_{d^n}$ be the iteration map with $d\geq 2$ and $n\geq 2$. There is a Zariski closed set $A\subset \Rat_d$,
which is  the preimage of the singularities of the variety
$\phi_{d,n}(\Rat_d)$, such that
$$\phi_{d,n}: \Rat_d\backslash A \rightarrow \Rat_{d^n}$$
is injective. Moreover, $A$ is a proper nonempty subset of $\Rat_d$.
\end{theorem}

\medskip

Finally, we characterize the condition that two non-exceptional rational functions share an iterate. Let $\Prep (f)=\{x\in \P^1~|f^n(x)=f^m(x), n>m\in \N\}$ be the set of preperiodic points of rational function $f$ and $\Per (f)=\{x\in \P^1~|f^n(x)=x, n\in \N^*\}$ be the set of periodic points of $f$, where $\N^*$ is the set of all positive integers.
\begin{theorem}\label{comthm}
Let $f$ and $g$ be non-exceptional rational functions  with degrees$\geq 2$. The following statements are equivalent:
\begin{itemize}
\item  $f$ and $g$ share an iterate, i.e. $f^n=g^m$ for some $n, m \in \N^*$.
\item There is some rational function $\varphi$ with degree$\geq 2$, such that $f\circ \varphi= \varphi \circ f$ and $g\circ \varphi=\varphi \circ g$.
\item $\mu_f=\mu_g$, and $J\cap \textup{Per}(f)\cap \textup{Per}(g)\neq \emptyset$.
\item $\Prep(f)=\Prep(g)$ and $J\cap \textup{Per}(f)\cap \textup{Per}(g)\neq \emptyset$.
\item $ \textup{Per}(f)=\textup{Per}(g)$.
\end{itemize}
\end{theorem}
The proof of Theorem \ref{comthm} uses the following results: for non-exceptional rational functions, Levin-Przytycki \cite{LP} showed that two rational functions having the same maximal measure should have the same set of preperiodic points. And conversely, Yuan and Zhang \cite{YZ} showed, via arithmetic methods, that rational functions having the same set of preperiodic points should have the same maximal measure.
\medskip

\noindent {\bf Historical background and related results.} In 1965, Brolin \cite{B4} introduced a probability measure $\mu_f$ for  polynomials $f$ with degree $\geq 2$. For any point $z_o\in \P^1$, with at most two exceptions, the sequence of sets $f^{-n}(z_o)$ equidistribute on the Julia set with respect to $\mu_f$, as $n\to \infty$. And for rational functions, this measure was introduced by Lyubich \cite{L3}  and independently  Freire-Lopez-Ma\={n}\'{e} \cite{FLM} in 1983. They showed that $\mu_f$ is the unique $f$-invariant measure supported on the Julia set, and  achieving the maximal entropy $\log d$.

As a general question, what are rational functions with the same maximal measure? For any non-exceptional polynomial  $f$, it is easy to read the symmetry group $\Sigma_{J_f}:=\{\sigma(z)=az+b|~\sigma (J_f)=J_f\}$ from the expression of $f$. After changing  coordinates, we can assume that $f$ is a monic and centered polynomial ($f(z)=z^d+az^{d-2}+\cdots$). So we can write $f(z)=z^lg(z^n)$ with $g(0)\neq 0$ and maximal possible $n$. Then, whenever $f$ is non-exceptional, we have $\Sigma_{J_f}=\{\sigma(z)=\zeta z|\zeta^n=1\}$. Then from (\ref{mpoly}), the expression of $M_f$ is clear for  non-exceptional  polynomials $f$.

For any  rational function $f$, let  $g\in M_f$ and $\sigma\in \Sigma_{\mu_f}$, we have $\sigma\circ g$ and $g\circ \sigma$ are both in $M_f$. So from Levin's result of the finiteness of the set $M_f\bigcap \Rat_n$, $\Sigma_{\mu_f}$ is a finite set if $f$ is not conjugate to $z^{\pm d}$. However, for rational functions $f$, it is still not known how to get the symmetry group $\Sigma_{\mu_f}$ or $ \Sigma_{J_f}$ (the subgroup of $PSL_2(\C)$ preserving  $J_f$) from the expression of $f$; see Levin's paper \cite{L1} \cite{L2}, and some other related results in \cite{DoM} and \cite{U1}.
And for  rational functions, in 1997, Levin and Przytycki's paper \cite{LP} has the following result:
\begin{theorem}[Levin-Przytycki \cite{LP}]\label{lp}
 Let $f$ and $g$ be two non-exceptional rational functions. The following two are equivalent:
 \begin{itemize}
 \item $\mu_f = \mu_g$;
 \item There exist iterates $F$ of $f$ and $G$ of $g$, integers $M$, $N\geq 1$, and locally defined branches of $G^{-1}\circ G$ and $F^{-1}\circ F$ such that
\begin{equation}\label{eqlp}
(G^{-1}\circ G)\circ  G^M= (F^{-1}\circ F)\circ  F^N.
\end{equation}
\end{itemize}
\end{theorem}
By analytic continuation, locally defined $G^{-1}\circ G$ and $F^{-1}\circ F$ can be extended to  multi-valued functions, acting by permuting the fibers of $G$ and $F$. Equation (\ref{eqlp}) implies that $f$ and $g$ have the same set of preperiodic points. Then as a consequence of Theorem 3 in Levin's paper \cite{L1} and Levin-Przytycki's theorem Theorem \ref{lp}, we have the following theorem:

\begin{theorem}[Levin-Przytycki]\label{main1}
 Let $f$ and $g$ be two non-exceptional rational functions with degrees$\geq 2$. Then  $\mu_f=\mu_g$ if and only if there are some iterates $F$ and $G$ of $f$ and $g$ such that
\begin{equation}\label{best}
F\circ F=F\circ G \textup{ and } G\circ F=G \circ G.
\end{equation}
\end{theorem}
Although the above theorem comes directly from Theorem 3 in \cite{L1} and Theorem \ref{lp}, we will  provide an easy proof later in Section 3 just by  using Levin-Przytycki's theorem (Theorem \ref{lp}).

So far as we know, (\ref{best}) is the strongest algebraic relation satisfied for all non-exceptional rational functions $f$ and $g$ with $\mu_f=\mu_g$.

\medskip
\noindent {\bf Outline of article.} The structure of this article goes as follows: in Section 2, we study the graph of the multi-valued functions $G^{-1}\circ G$, and give three theorems Theorem \ref{genus0}, \ref{3333} and \ref{111}, which will be used in proving Theorem \ref{generic rational} and \ref{mainthm} in later sections.  In Section 3 and 4, we prove the main theorems stated in this section.  Finally, in Section 5, we will discuss  rational functions with common iterates and prove Theorem \ref{comthm}.

In the writing of this paper, we learned of two related articles in preparation. Related to Theorem \ref{injective}, Adam Epstein has shown that the iteration map $\phi_{d,n}$ ia an immersion for all $d$ and $n\geq 2$; see Proposition \ref{singular}. Related to Theorem \ref{counter}, when $T$ is assumed to be a polynomial, Avanzi, Zannier, Carney, Hortsch and Zieve has provided a complete list of such examples; See \cite{AZ} and \cite{CHZ}.

\bigskip

\noindent {\bf Acknowledgements.} I would like to thank Chong Gyu Lee, Adam Epstein,  Michael Zieve and especially Laura DeMarco for lots of helpful comments and suggestions.

%%%%%%
%%%%%%

\section{Graph of the multi-valued functions $G^{-1}\circ G$}\label{sec2}
In this section, we study the geometry of the algebraic curves defined by $$V_G=\{(x,y)\in \P^1\times\P^1| G(x)=G(y)\}\subset \P^1\times\P^1$$ for rational functions $G$. The irreducible components of $V_G$ correspond to multi-valued functions $G^{-1}\circ G$, as appearing in equation (\ref{eqlp}). We prove Theorem \ref{genus0} , \ref{3333} and \ref{111} allowing us to estimate the genus of the irreducible components of $V$.

\subsection{ Multi-valued functions $G^{-1}\circ G$} For a rational function $G$ with degree $d\geq 2$, we consider the correspondence  $V_G=\{(x,y)\in \P^1\times\P^1| G(x)=G(y)\}\subset \P^1\times\P^1$ of function $G$. Obviously, $V_G$ is a projective variety of bidegree ($d,d$), consisting of finitely many irreducible curves.
Let $V_o$ be an irreducible component of $V_G$ with bidegree ($r_1, r_2$). Geometrically,  for $i=1$ and $ 2$, $r_{2-i}$ is the topological degree of the coordinate projection map $\pi_i: V_o\to \P^1$. The diagonal $\bigtriangleup \subset \P^1\times\P^1$ is an irreducible
component of $V_G$ with bidegree ($1,1$).

For any two distinct and noncritical points $z_1,z_2$ with $G(z_1)=G(z_2)$, there is a unique  holomorphic germ  $b(z)$ locally defined at $z_1$, such that $b(z_1)=z_2$ and $G\circ b=G$. By analytic continuity, we can extend this germ $b(z)$ to  a multi-valued function from the function on all of $\P^1$, denoted as $G^{-1}\circ G$. To be clear, throughout this paper,  $G^{-1}\circ G$ will refer to a particular  multi-valued functions defined in this way. Any such multi-valued function $G^{-1}\circ G$ corresponds to an irreducible component $V_o$ of  $V_G$. Conversely, each irreducible component $V_o$ corresponds to exactly one multi-valued function $G^{-1}\circ G$. We call $V_o$ the {\bf graph} of its corresponding multi-valued function $G^{-1}\circ G$. And then $V_G$ is the union of the graphes for all the multi-valued functions $G^{-1}\circ G$.

Let $\Crit (G)\subset \P^1$ be the set of  critical points of $G$ and $\widetilde{\Crit}(G)$ be the preimage of the set of critical values of $G$. Let $S_G=\P^1 \setminus \widetilde{\Crit}(G)$ and $\widetilde{S}_G$ be the universal cover of $S_G$, i.e. we have the covering map
\begin{equation}\label{group}
\widetilde{S}_G \longrightarrow  S_G=\widetilde{S}_G/H,
\end{equation}
where $H$ is a subgroup of the automorphism group of $\widetilde{S}_G$, and $H$ is isomorphic to the fundamental group of $S_G$.

Fix a non-diagonal irreducible component $V_o$ of $V_G$ and its corresponding multi-valued function $G^{-1}\circ G$. Let ($r_1, r_2$) be the bidegree of  $V_o$ in $\P^1\times\P^1$. $G^{-1}\circ G$ is a multi-valued function from $S_G$ to $S_G$. Although it may not be single-valued, we can lift it to the universal cover, and get a single valued function $h_o$ from $\widetilde{S}_G$ to $\widetilde{S}_G$.  $h_o$ is an automorphism of
$\widetilde{S}_G$ without fixed point, since $V_o$ is not the diagonal. Moreover, we have $h_o \notin H$, otherwise $G^{-1}\circ G$ would be identity map, here
  $H$ is the group in (\ref{group}). Now, we can use the index of the  fundamental groups to interpret $r_1$ and $r_2$. For any $\widetilde{x}\in \widetilde{S}_G$, $H\widetilde{x}\in S_G=\widetilde{S}_G/H$. For any $h_i, h_j \in H$,
  $Hh_oh_i\widetilde{x}=Hh_oh_j\widetilde{x}$ if and only if $Hh_oh_i=Hh_oh_j$, if and only
   if $h_ih_j^{-1}\in h_o^{-1}Hh_o$. Since $V_o$ is of  bidegree ($r_1, r_2$), each coset $H\widetilde{x}$ splits into $\{H_i\widetilde{x}\}$, where $H_i$ are the cosets of $H\cap
   h_o^{-1}Hh_o$ in $H$, for $i=1,2, \cdots, r_2$.  As a consequence, we can write
   $r_2=[H: H\cap h_o^{-1}Hh_o]$. In order to find $r_1$, we can use $h_o^{-1}$ instead of $h_o$. Similarly, we have  $r_1=[H: H\cap
   h_oHh_o^{-1}]$. The map $G: S_G \to \P^1\setminus CV(G)$ is a covering map of degree $d$, where $d$ is the degree of $G$ and $CV(G)$ is the set of critical values of $G$.
$$\widetilde{S}_G\longrightarrow S_G=\widetilde{S}_G/H \longrightarrow \P^1\setminus CV(G)=\widetilde{S}_G/\widetilde{H}$$
Because $G^{-1}\circ G$ permutes the points in each fiber of $G$, we have $h_o\in \widetilde{H}$ and $[\widetilde{H}:H]=d$.
$$
\begin{array}{lll}
r_1d
&=[\widetilde{H}:H\cap h_oHh_o^{-1}]\\[6pt]
&=[h_o^{-1}\widetilde{H}h_0:h_o^{-1}(H\cap h_oHh_o^{-1})h_o]\\[6pt]
&=[\widetilde{H}:h_o^{-1}Hh_o\cap H]\\[6pt]
&=r_2d
\end{array}
$$
Then we have $r_1=r_2$.

\begin{prop}
Let $G$ be a rational function with degree $d\geq 2$, $V_G=\{(x,y)\in \P^1 \times \P^1 | G(x)=G(y)\}$. The graph $V_G$ is of bidegree $(d,d)$, and any irreducible component (as a variety) of $V_G$ is of bidegree $(r,r)$ with $1\leq r\leq d-1$. Moreover, the sum of  bidegrees $r$ of the irreducible components of $V_G$ is $d$.

\end{prop}

Let $V_o$ be an irreducible component of $V_G$ with bidegree ($r,r$). It may contain singularities, but we can normalize it and get a smooth curve $\widetilde{V}_o$ as its normalization. Then we have the following natural projections:
$$\xymatrix{ \widetilde{V}_o  \ar[r]^{\pi} &V_o \ \ar[r]^{\pi_i} & \P^1},$$
where, $\pi_i$ is the coordinate projection, for $i=1,2$.

We use $d_{h,x}$ to denote  the local degree of a holomorphic map $h(z)$ at point $x$.   For any $(\widetilde{x},\widetilde{y})\in \widetilde{V}_o$ with $\pi((\widetilde{x},\widetilde{y}))=(x,y)\in V_o$, we can express the local degree of the map $\pi_1\circ \pi$ at $(\widetilde{x},\widetilde{y})$ in terms of $d_{G,x}$ and $d_{G,y}$. From the local behavior of $G$ at points $x$ and $y$, it has
\begin{equation}\label{degree}
d_{\pi_1\circ \pi, (\widetilde{x},\widetilde{y})}=\frac{d_{G,y}}{\textup{gcd}(d_{G,x},d_{G,y})}.
\end{equation}

\subsection{Genus zero components of the graph $V_G$}
\begin{theorem}\label{genus0}
Let $G^{-1}\circ G$  be a multi-valued function with corresponding irreducible component $V_o\subset V_G$. The following two are equivalent:
\begin{itemize}
\item The normalization $\widetilde{V}_o$ of $V_o$ has  genus  zero.
\item There exist  rational functions $\widetilde{G}$ and $\widetilde{F}$ such that $(G^{-1}\circ G) \circ \widetilde{G}=\widetilde{F}$.
\end{itemize}
\end{theorem}
Proof: Assume that there are rational functions $\widetilde{G}$ and $\widetilde{F}$, such that $(G^{-1}\circ G) \circ \widetilde{G}=\widetilde{F}$. Then we have a well defined map
$$\rho: \P^1 \to V_o$$
with $\rho(z)=(\widetilde{G}(z), \widetilde{F}(z))$. The map $\rho$ can be lifted to $V_o$'s normalization $\widetilde{V}_o$, and denote the lifting map as $\widetilde{\rho}$,
$$\widetilde{\rho}: \P^1 \to \widetilde{V}_o.$$
 The lifting map $\widetilde{\rho}$  is  holomorphic from $\P^1$ to the smooth curve $\widetilde{V}_o$. And by Riemann-Hurwitz formula, there is no nonconstant holomorphic map from $\P^1$ to a curve with genus greater than zero. Then the genus of $\widetilde{V}_o$ should be zero.

\medskip
Conversely, if the genus of the smooth curve $\widetilde{V}_o$ is zero.  We can parameterize $\widetilde{V}_o$ by $\P^1$ using some parametrization $\widetilde{\rho}$,
$$\widetilde{\rho}: \P^1 \to \widetilde{V}_o.$$
After projecting it down to $V_o$, we get the following map:
$$\rho=\widetilde{\rho}\circ \pi: \P^1 \to V_o,$$
where $\rho(z)=(\widetilde{G}(z), \widetilde{F}(z))$ with $\widetilde{G}$ and $ \widetilde{F}$ being  rational functions of degree $r$. From this parametrization, it is clear that
$$(G^{-1}\circ G) \circ \widetilde{G}=\widetilde{F}.$$\qed

{\bf Example.} Let $T(z)=z^3-3z$ be a degree $3$ polynomial. Easy to check that the graph $V_T$ has two irreducible components. One of them is the diagonal of $\P^1\times \P^1$ with bidegree $(1,1)$ corresponding to the identity map, and the other one $V_o$ is of bidegree $(2,2)$. From Riemann-Hurwitz formula, it can be computed that the genus of $V_o$'s normalization $\widetilde{V}_o$ is zero. \qed

\subsection{Functions $G$ without nontrivial genus zero components of  $V_G$}\begin{theorem}\label{3333}
Let $G$ be a rational function with degree $d_G\geq 3$. Assume that there are at least three simple critical values for $G$. Then for any irreducible component $V_o\subset V_G$ with
 bidegree $(r,r\geq 2)$, its normalization $\widetilde{V}_o$ has  genus$\geq 1$.
\end{theorem}
\proof By assumption, $V_o\subset V_G$ is an irreducible component of $V_G$ with bidegree ($r,r\geq 2$), corresponding to some multi-valued function $G^{-1}\circ G$.
  Let $\{y_1,y_2,y_3\}$ be  three simple critical values of the rational function $G$. Consider the sets:
$$\{x_{i,1}, x_{i,2},\cdots, x_{i,d_G-1}\}=G^{-1}(y_i), $$
 where $x_{i,d_G-1}$ is a critical point of $G$ and $i=1,2,3$. Since $V_o$ has bidegree ($r,r$),  there are at least $r-1$ noncritical points in $G^{-1}(y_i)$, which can be assumed to be $\{x_{i,1}, x_{i,2}, \cdots, x_{i,r-1}\}$, such
that
$$(x_{i,j},x_{i,d_G-1})\in V_o, \textup{ for } j=1,2,\cdots, r-1.$$
For each of such points $(x_{i,j},x_{i,d_G-1})\in V_o$ with $j\leq r-1$, since $d_{G, x_{i,d_G-1}}>1$ and  $d_{G, x_{i,j}}=1$ ,by equation (\ref{degree}), $d_{\pi_1\circ \pi, (\widetilde{x}_{i,j},\widetilde{x}_{i,d_G-1})}>1$ for any $(\widetilde{x}_{i,j},\widetilde{x}_{i,d_G-1})$ in $\pi^{-1}((x_{i,j},x_{i,d_G-1}))\subset
\widetilde{V}_o$. Hence $(\widetilde{x}_{i,j},\widetilde{x}_{i,d_G-1})$ is a critical point for the projection map $\pi_1\circ \pi$:
$$\xymatrix{ \widetilde{V}_o  \ar[r]^{\pi} &V_o \ \ar[r]^{\pi_1} & \P^1}.$$
So we have at least $3*(r-1)$ critical points for the map $\pi_1\circ \pi$. By Riemann-Hurwitz formula, we have
$$
\begin{array}{lll}
2-2\textup{genus}(\widetilde{V}_o)
&=2r-\#\textup{ of critical points of $\widetilde{V}_o$}\\[6pt]
&\leq 2r-3*(r-1)\\[6pt]
&\leq 3-r.
\end{array}
$$
Since $r\geq 2$, we have
$$\textup{genus}(\widetilde{V}_o)\geq (r-1)/2> 0.$$ \qed

For degree $2$ case, we have the following:
\begin{theorem}\label{111}
Let $g$ be a degree two rational function with two disjointed critical orbits and none of the critical points is preperiodic.  Let $G=g^n$ for some $n\geq 1$.
Assume $V_o$   is an irreducible component of $V_G$ such that the corresponding $G^{-1}\circ G$ is not the identity map or $\sigma_f$. Then its normalization $\widetilde{V}_o$ has  genus$\geq 1$.
\end{theorem}

\proof Let ($r,r$) be the bidegree of  $V_o\subset V_G$. Let $x_{1,1}, x_{2,1}$ be the two critical points of $f$ and
$$\{x_{i,1}, x_{i,2},\cdots, x_{i,d_G-1}\}=G^{-1}(G(x_{i,1})), \textup{ for $i=1,2.$}$$
 As the multi-valued function $G^{-1}\circ G$ does not correspond to the identity map or $\sigma_f$, there are $r$ noncritical points in $\{x_{i,2},\cdots, x_{i,d_G-1}\}$, which can be assumed to be $\{x_{i,2}, x_{i,3}, \cdots, x_{i,r+1}\}$,
 such
that
$$(x_{i,j},x_{i,1})\in V_o, \textup{ for } j=2,3,\cdots, r+1.$$
For each of such points $(x_{i,j},x_{i,1})\in V_o$ with $1<j\leq r+1$, since $d_{G, x_{i,1}}>1$ and  $d_{G, x_{i,j}}=1$ ,by equation (\ref{degree}), $d_{\pi_1\circ \pi, (\widetilde{x}_{i,j},\widetilde{x}_{i,1})}>1$ for any $(\widetilde{x}_{i,j},\widetilde{x}_{i,1})$ in $\pi^{-1}((x_{i,j},x_{i,1}))\subset
\widetilde{V}_o$. Hence $(\widetilde{x}_{i,j},\widetilde{x}_{i,1})$ is a critical point for the projection map $\pi_1\circ \pi$:
$$\xymatrix{ \widetilde{V}_o  \ar[r]^{\pi} &V_o \ \ar[r]^{\pi_1} & \P^1}.$$
So we have at least $2r$ critical points for the map $\pi_1\circ \pi$. By Riemann-Hurwitz formula, we have
$$
\begin{array}{lll}
2-2\textup{genus}(\widetilde{V}_o)
&=2r-\#\textup{ of critical points of $\widetilde{V}_o$}\\[6pt]
&\leq 2r-2r\\[6pt]
&=0.
\end{array}
$$
Then we have
$$\textup{genus}(\widetilde{V}_o)\geq 1.$$ \qed

\section{Rational functions with common measure of maximal entropy}
In this section, we study the relation of two rational functions $f$ and $g$ with $\mu_f=\mu_g$, and then prove Theorem \ref{mainthm}, \ref{counter} and \ref{main1}. Moreover, we give examples of non-exceptional functions  $f$ and $g$ with $\mu_f=\mu_g$ and they do not satisfy (\ref{polycom}).

\begin{theorem}\label{2222}
Assume $f$ and $g$ are rational functions with degrees$\geq 2$,  satisfying
\begin{equation}\label{lpequation}
(g^{-1}\circ g) \circ g=(f^{-1}\circ f)\circ f
\end{equation}
for some multi-valued functions $g^{-1}\circ g$ and $f^{-1}\circ f$. Then there are some iterates $F$ and $G$ of $f$ and $g$, such that:
$$F\circ F=F\circ G, ~G\circ F=G\circ G.$$
\end{theorem}
\proof Choose a point $a_0\in \P^1$ such that $a_0$ is neither in  $f$'s critical orbits nor in $g$'s critical orbits. Let $a_0, a_1, a_2, \cdots$ be a sequence of points such that $g(a_i)=a_{i-1}$. From equation (\ref{lpequation}), for any $i\geq 1$, after composing each side $i$ times with themself, we have
\begin{equation}\label{cc}
(g^{-1}\circ g) \circ g^i=(f^{-1}\circ f)\circ f^i \textup{ or } f\circ (g^{-1}\circ g) \circ g^i= f \circ f^i.
\end{equation}
Then for each $i\geq 1$, there is function germ $(g^{-1}\circ g)_i$ of $g^{-1}\circ g$, locally defined at $a_0$, such that functions germs
$$f\circ (g^{-1}\circ g)_i\circ g^i|_{\textup{near $a_i$}}=f\circ f^i|_{\textup{near $a_i$}},$$
which are locally defined near $a_i$. Let $b_i=(g^{-1}\circ g)_i(a_0)$. Since $g(b_i)= g(a_0)$ for any $i\geq 1$, there are only finitely many distinct $b_i$.
Choose some $j>2i_1>2$ such that $b_{i_1}=b_j$. Then we have germs $(g^{-1}\circ g)_{i_1}=(g^{-1}\circ g)_j$ locally defined near $a_0$. From equation (\ref{cc}), we have locally defined germs
$$f\circ (g^{-1}\circ g)_{i_1}\circ g^{i_1}|_{\textup{near $a_{i_1}$}}=f\circ f^{i_1}|_{\textup{near $a_{i_1}$}}$$
$$f\circ (g^{-1}\circ g)_j\circ g^j|_{\textup{near $a_j$}}=f\circ f^j|_{\textup{near $a_j$}}$$
As $(g^{-1}\circ g)_{i_1}=(g^{-1}\circ g)_j$, combining above two equations, it follows
$$f\circ f^{i_1}\circ g^{j-i_1}|_{\textup{near $a_j$}}=f\circ (g^{-1}\circ g)_{i_1}\circ g^{i_1}\circ g^{j-i_1}|_{\textup{near $a_j$}}=f\circ f^j|_{\textup{near $a_j$}}$$
Because both sides of the above equation are germs of rational functions, we have
$$f\circ f^{i_1}\circ g^{j-i_1}=f\circ f^j.$$
Since $j\geq 2i_1\geq 2$, we can post compose some iterate of $f$ to both sides of the above equation:
$$ f^{j-i_1}\circ g^{j-i_1}=f^{j-i_1}\circ f^{j-i_1}.$$
Let $f_o=f^{j-i_1}$ and $g_o=g^{j-i_1}$.  The above equation shows $(f_o^{-1}\circ f_o)\circ f_o= g_o$. Consequently, $(f_o^{-1}\circ f_o)\circ f_o^i=g_o^i$ and $f_o^i\circ f_o^i=f_o^{i}\circ g_o^i$ for any $i\geq 1$.

Since we have $(f_o^{-1}\circ f_o)\circ f_o= g_o$, repeating the same process of the above proof, there is some $i_o\geq 1$ such that $g_o^{i_o}\circ f_o^{i_o}=g_o^{i_o}\circ g_o^{i_o}.$

Let $F=f_o^{i_o}$ and $G=g_o^{i_o}$. Then we have $F\circ F=F\circ G, ~G\circ F=G\circ G$.\qed

\medskip
As a consequence of Theorem \ref{poly id}, we can easily prove Theorem \ref{main1}   by just using Theorem \ref{lp}.

{\bf Proof of Theorem \ref{main1}.} From Theorem \ref{lp}, there are some iterates $f_o$ and $g_o$ of $f$ and $g$ and $M, N\geq 1$, such that
$$(g_o^{-1}\circ g_o)\circ g_o^M=(f_o^{-1}\circ f_o)\circ f_o^N.$$
Since for any multi-valued function $(g_o^{-1}\circ g_o)$, we can choose a multi-valued function $(g_o^{-M}\circ g_o^M)$ which the same as $(g_o^{-1}\circ g_o)$, from above equation, we have
$$(g_o^{-M}\circ g_o^M)\circ g_o^M=(f_o^{-N}\circ f_o^N)\circ f_o^N.$$
By Theorem \ref{2222}, there are iterates $F$ and $G$ of $f_o^N$ and $g_o^M$, such that
$$F\circ F=F\circ G, ~G\circ F=G\circ G.$$

\qed

Given any non-exceptional polynomial $g$ with degree$\geq 2$, for any rational function $f$ with $\mu_g=\mu_f$, it was known that $f$ should also be a polynomial; see \cite{OS}. As a corollary of Theorem \ref{main1}, here we give an easy proof of this result.

\begin{cor}\label{poly id}
Let $g$ be a non-exceptional  polynomial with degree $d\geq 2$. Then any rational function $f$  with $\mu_f=\mu_g$ should be a polynomial. Consequently, there exist some $m,n\geq 1$, s.t.
$$f^n=\sigma \circ g^m,$$
where $\sigma(z)=az+b$ is an affine transformation preserving $\mu_g=\mu_f$.
\end{cor}
\proof From Theorem \ref{main1}, there are some iterates  $F$ and $G$ of $f$ and $g$, such that
$$G\circ F=G\circ G.$$
Exception set of a rational function $h$ is the  maximal  finite set, which is invariant under $h$. Exception set can only be an empty, one point or two points set. If the exception set is one point, then $h$ is
conjugate to a polynomial. Since $G$ is a non-exceptional polynomial, its exception set is $\{\infty\}$. Then $\{\infty\}$ is also an invariant set of $F$, which means that $F$ is a polynomial. If the exceptional set of $F$ contains two points, then $F$ is conjugate to polynomial $z^{d_F}$, which means $F$ is an exceptional polynomial. And because $\mu_F=\mu_G$, $G$ is exceptional. This contradicts to the assumption. Consequently, $\{\infty\}$ is the exceptional set of $F$. Since $F$ is some iterate of $f$, they should have the same exceptional set. Then $f$ should be a polynomial. The last statement comes from the main theorem of \cite{SS}.\qed
\medskip

{\bf Proof of Theorem \ref{mainthm}. }  Because any exceptional rational function has at most two simple critical values, $f$ is non-exceptional. As $\mu_g=\mu_f$, then by Theorem \ref{main1}, there are some integers $m,n\geq 1$, such that for $F=f^n$ and $ G=g^m$,
\begin{equation}\label{bbbbb}
F\circ F=F\circ G, \textup{ or } (F^{-1}\circ F)\circ F=F\circ G
\end{equation}

If $F=G$,  $f$ and $g$ has a common iterate. Then the statement is satisfied. So we can assume that $F\neq G$.   Let $k\geq 0$ be the smallest integer such that $f^k\circ F\neq f^k\circ G$ and $f^{k+1}\circ F=f^{k+1}\circ G$. Since
\begin{equation}\label{bbbb}
(f^{-1}\circ f)\circ f^k\circ F= f^k\circ G,
\end{equation}
by Theorem \ref{genus0} and \ref{3333},  the corresponding irreducible component of the multi-valued function $(f^{-1}\circ f)$ in (\ref{bbbb}) should have bidegree ($1,1$). It means the multi-valued function $\sigma=(f^{-1}\circ f)$ is a M$\ddot{\textup{o}}$bius transformation and $f\circ \sigma=f$.

Under changing of coordinates, we can assume that $\sigma(z)=\zeta z$ where $\zeta$ is a $k$'s primitive root of unit. If $k\geq 2$, then we can decompose $f(z)$ into
$$f(z)=f_o(z^k).$$
Since $k \geq 2$, from the above decomposition, $f$ cannot have three simple critical values. This is a contradiction. So we have
$$\sigma(z)=(f^{-1}\circ f)(z)=z.$$
And by (\ref{bbbb}), finally we get $f^k\circ F=f^k\circ G$, which is a contradiction to the assumption that $f^k\circ F\neq f^k\circ G$. In all, it has $F=G$, i.e. $f$ and $g$ share an iterate. \qed

\medskip
{\bf Remark.} Theorem \ref{mainthm} asserts that for general $f\in \Rat_d$ with degree $d\geq3$, $\mu_f=\mu_g$ implies that $f$ and $g$ share an
iterate. And as we discussed in the introduction, the existence of the special symmetry $\sigma_f$ for any $f\in \Rat_2$ prevents the same conclusion
as in Theorem \ref{mainthm}. Precisely, for any $f\in \Rat_2$ and $g=\sigma_f\circ f$, we have  $\mu_f=\mu_g$,  but they never share an iterate. However, we can modify it a bit, and show that for generic $f\in \Rat_2$ (see Theorem \ref{111}), $\mu_f=\mu_g$ implies that
$g^m=\sigma_f\circ f^n$ or $f^n$; for details see the proof of Theorem \ref{generic rational}.

\medskip

However, $\mu_f=\mu_g$ does not always imply that $f$ and $g$ share an iterate. Even worse, Theorem \ref{counter} asserts that $f$ and $g$ may not even satisfy (\ref{polycom}) for any M$\ddot{\textup{o}}$bius transformation $\sigma$.

{\bf Proof Theorem \ref{counter}.} Since $T\circ R=T\circ S$,  we have $f\circ f=f\circ g$ and $T\circ f^i(z)=T\circ g^i(z)$ for any $i\geq 1$. Consequently,
 $$(f^{-1}\circ f)\circ f=g.$$ Then from Theorem \ref{lp}, $\mu_f=\mu_g$.

Assume that there exist integers $n,m\geq 1$ and $\sigma\in PSL_2(\C)$ such that
$$f^n=\sigma\circ g^m.$$
Since $f$ and $g$ have the same degree,  we have $n=m$, i.e. $f^n=\sigma\circ g^n.$
$$f^n(z)=R\circ T\circ f^{n-1}=R\circ T\circ g^{n-1}=\sigma\circ g^n=\sigma \circ S\circ T\circ g^{n-1}$$
So we have $R=\sigma \circ S$, which contradicts to the assumption in this theorem.

For the first statement of this theorem, see the following example.\qed
\medskip

\begin{figure}
\begin{minipage}[t]{0.5\linewidth}\centering
\includegraphics[width=3in]{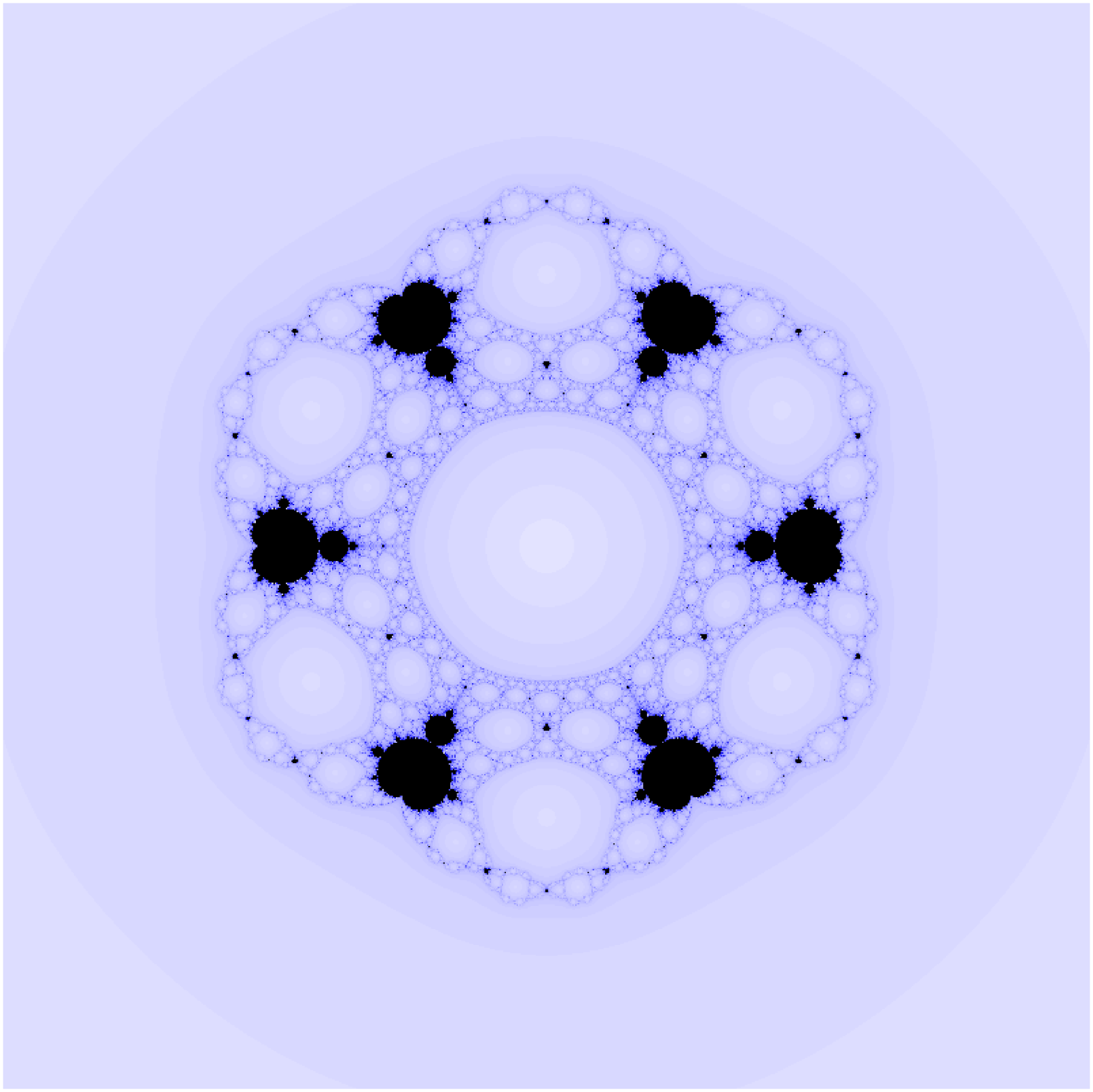}\caption{ The parameter space of $f_a$}\label{fig1}\end{minipage}
\begin{minipage}[t]{0.5\linewidth} \centering
\includegraphics[width=3in]{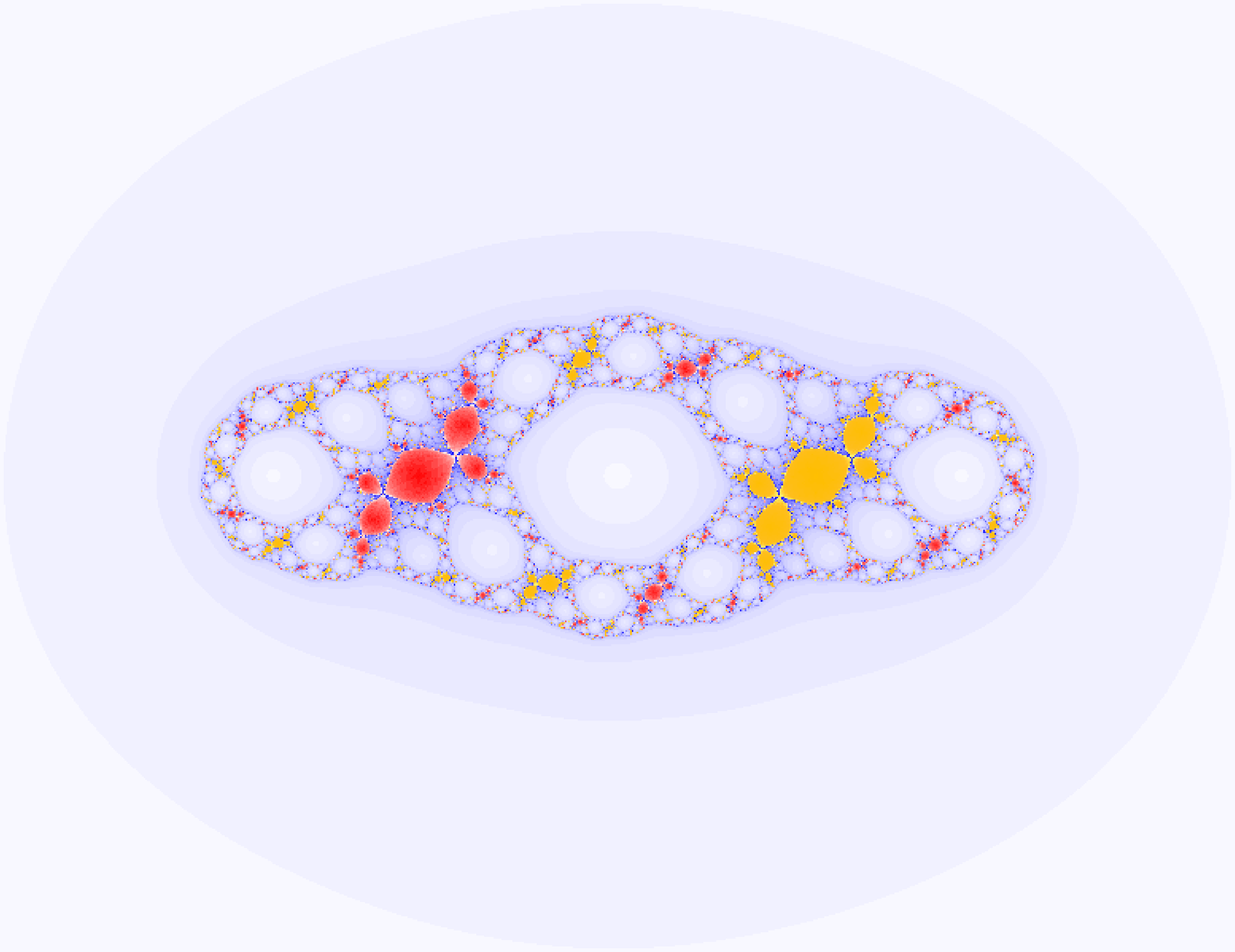}
\caption{Julia set of $f_{a}$ with $a=0.4843+0.07776i$}
\label{fig2}
\end{minipage}
\end{figure}
\medskip

{\bf Example.} To illustrate Theorem \ref{counter}, let $T(z)=z^3-3z, ~R(z)=az+\frac{1}{az}$ and $S(z)=a\omega z+\frac{1}{a\omega z}$, with $\omega^2+\omega +1=0$ and $a\in \C^*$. For any $a\in \C^*$, it easy to check that $T\circ R=T\circ S$ and there is no $\sigma \in PSL_2(\C)$ such that $R=\sigma\circ S$. So from Theorem \ref{counter}, we know that $f=R\circ T$ and $g=S\circ T$ have the same measure of maximal entropy.  And for any $n,m\geq 1, \sigma \in PSL_2(\C)$, we have
$$f^n\neq \sigma \circ g^m.$$
It is not hard to see that neither $f$ nor $g$ is exceptional rational function, since they are not critical finite. These are not the only examples satisfying the assumptions in Theorem \ref{counter}, for example, Michael Zieve suggests the following functions
   $$T(z)=z^n(z+1)^m, ~ R(z)=(1-z^n)/(z^{n+m}-1), ~ S(z)=z^m(1-z^n)/(z^{n+m}-1).$$
 There are more such $T, R$ and $S$ satisfying the assumptions in Theorem \ref{counter}; see \cite{AZ} and \cite{CHZ}.\qed

{\bf Remark.} In the above example,  rational functions $f_a(z)=a(z^3-3z)+\frac{1}{a(z^3-3z)}$ and $g_a(z)=f_{a\omega}(z)$ come from composition of rational functions $T(z)=z^3-3z, R(z)=az+\frac{1}{az}$ and $S(z)=a\omega z+\frac{1}{a\omega z}$.  And $T, R, R$ satisfy the assumptions in  Theorem \ref{counter}. Figure \ref{fig1} is the  parameter space of $f_a$ which indicates that $\mu_{f_{\zeta a}}=\mu_{f_a}$ for any $\zeta$ with $\zeta^6=1$. Actually, by Theorem \ref{counter}, we know that $\mu_{f_{\omega a}}=\mu_{f_a}$ for $\omega^2+\omega+1=0$ and $f_a^2=f^2_{-a}$. So $\mu_{f_{\zeta a}}=\mu_{f_a}$ for any $\zeta$ with $\zeta^6=1$. Since $\infty$ is a supper attracting point and $f_a$ is not a polynomial, there is a critical point attracted to $\infty$ and it is not periodic. As exceptional functions are all post-critical finite, $f_a$ won't be exceptional. By Theorem \ref{counter}, $f_a$ and $g_a$ has the same measure and there is no iterate of $f_a$ conjugated to an iterate of $g_a$. However, for any non-exceptional polynomials $f$ and $g$ with $\mu_f=\mu_g$, there always exist iterates of $f$ and $g$ which are in the same conjugacy class.

There are more examples of such rational functions $T, R$ and $S$ as in Theorem \ref{counter}. For example, let $t, r, s$ be rational functions satisfying assumptions in Theorem \ref{counter}. Then for any rational function $h$,  $T=h\circ t, R=r, S=s$ satisfy the same assumptions. It has been classified for all such rational functions $T, R$ and $S$, with the restriction that $T$ is a polynomial. However, when $T$ is not a polynomial, it is still not known how to classify it. For details, please refer to \cite{AZ} and \cite{CHZ}.

\section{Generic rational function with identical measure}

In this section, we are going to prove Theorem \ref{injective}, which indicates iteration map is  one-to-one for general points. And then get to prove the main theorem, Theorem \ref{generic rational}, which says: for generic  rational functions $f\in \Rat_d$,
 we have  $M_f=\{f, f^2, f^3, \cdots\}$ or $ \{f, \sigma_f\circ f, f^2, \sigma_f\circ f^2, \cdots\}$, where $M_f$ the set of rational
 functions with the same maximal entropy measure as $f$.

\medskip
For $d, n\geq 2$, let $x$ be a point in $ \Rat_d$. There is an induced map between the tangent spaces of $x\in \Rat_d$ and $\phi_{d,n}(x)\in \Rat_{d,n}$:
$$\phi_{d,n~*}: T_x \to T_{\phi_{d,n}(x)}.$$
The map $\phi_{d,n}$ is singular at $x\in \Rat_d$ if the induced map between the tangent spaces $T_x$ and $T_{\phi_{d,n}(x)}$
 is not injective. The map $\phi_{d,n}$ is nonsingular if it is not singular at any point of $\Rat_d$.
\medskip
For any $f\in \Rat_d$, we can express it as
$$f(z)=\frac{h(z)}{k(z)}=\frac{a_dz^d+a_{d-1}z^{d-1}+\cdots+a_0}{b_dz^d+b_{d-1}z^{d-1}+\cdots+b_0}.$$
We can define a holomorphic map from $t \in \C$ to $f_t\in \Rat_d$ in a neighborhood of $0\in \C$:
$$f_t(z)=\frac{(a_d+t\alpha_d)z^d+(a_{d-1}+t\alpha_{d-1})z^{d-1}+\cdots+(a_0+t\alpha_0)}{(b_d+t\beta_d)z^d+(b_{d-1}+t\beta_{d-1})z^{d-1}+\cdots+
(b_0+t\beta_0)}.$$
It can be checked that this parametrization map is singular at $t=0$ iff
$$(\alpha_dz^d+\alpha_{d-1}z^{d-1}+\cdots+\alpha_0)k(z)-(\beta_dz^d+\beta_{d-1}z^{d-1}+\cdots+\beta_0)h(z)=0.$$
Take the derivative of $f_t$ with respect to $t$,
$$\frac{d f_t(z)}{dt}|_{t=0}= \frac{(\alpha_dz^d+\alpha_{d-1}z^{d-1}+\cdots+\alpha_0)k(z)-(\beta_dz^d+\beta_{d-1}z^{d-1}+\cdots+\beta_0)h(z)}
{(b_dz^d+b_{d-1}z^{d-1}+\cdots+b_0)^2}.$$
From the above expression, the map $t\to f_t$  is singular at $0$ iff $\frac{d f_t(z)}{dt}|_{t=0}=0$.

\medskip
More generally, for any holomorphic map ($t\to f_t$) from a neighborhood of $0\in \C$ to $\Rat_d$, we can express it as
$$f_t(z)=\frac{a_d(t)z^d+a_{d-1}(t)z^{d-1}+\cdots+a_0(t)}{b_d(t)z^d+b_{d-1}(t)z^{d-1}+\cdots+b_0(t)},$$
where $a_i(t)$ and $b_i(t)$ are holomorphic. Similar to the special case we discussed above, it can be checked that
\begin{equation}\label{iffsingular}
 f_t\textup{ is singular at $t=0$}  \Longleftrightarrow \frac{d f_t(z)}{dt}|_{t=0}=0.
\end{equation}
When the map  $t \to f_t$ is nonsingular at $t=0$, then $\frac{d f_t(z)}{dt}|_{t=0}$ is
a nonzero rational function with degree at most $2d$.

The following proposition has been proved by Adam Epstein \cite{Ep}. For completeness, we will prove it again here.

\begin{prop}\label{singular}
The map $\phi_{d,n}: \Rat_d \to \Rat_{d^n}$ is nonsingular. In particular, $\phi_{d,n}$ is an immersion from $\Rat_d$ to $\Rat_{d^n}$.
\end{prop}
\proof In order to prove that $\phi_{d,n}$ is nonsingular, it  suffices to prove that if a holomorphic map $t\to f_t$ from
a neighborhood of $0\in \C$ to $\Rat_d$ is nonsingular at $t=0$, then the map $t \to f^n_t$ is nonsingular at $t=0$.

First assume that the holomorphic map $t\to f_t$ is singular at $t=0$. Let $z_0\in \P^1$ be a periodic point of $f_0$ with period $p\geq 1$ and
multiplier $\frac{df^p_0}{dz}(z_0)\neq 1$. Then there is a holomorphic motion $z_t$ of the periodic point $z_0$ in $\P^1$ such that
\begin{equation}\label{pperiodic}
f_t^p(z_t)=z_t.
\end{equation}
We claim that $\frac{dz_t}{dt}|_{t=0}=0$, i.e. the holomorphic motion $z_t$ of the periodic point $z_0$ is singular at $t=0$. Indeed, let
$$\psi_p(z_1, z_2)=f_{z_2}^p(z_1),$$
and then $\psi_p(z,t)=f_{t}^p(z)$. By taking the derivative of $t$ for both sides of the equation (\ref{pperiodic}),
\begin{equation}\label{partial}
\frac{d\psi_p(z_t, t)}{dt}=\frac{\partial \psi_p}{\partial z_1}\frac{\partial z_t}{\partial t}+\frac{\partial\psi_p}{\partial z_2}=\frac{dz_t}{dt}
\end{equation}
Since $t\to f_t$ is singular at $t=0$, the map $t\to f_t^p$ is singular at $t=0$. As a consequence of (\ref{iffsingular}), $\frac{\partial\psi_p}{\partial z_2}$ is zero at
$(z_0, 0)$. Then from equation (\ref{partial}):
$$\frac{df^p_0}{dz}(z_0)\frac{\partial z_t}{\partial t}|_{t=0}=\frac{\partial \psi_p}{\partial z_1}(z_0,0)\frac{\partial z_t}{\partial t}|_{t=0}=
\frac{dz_t}{dt}|_{t=0}.$$
And by the assumption that $\frac{df^p_0}{dz}(z_0)\neq 1$, we have $\frac{dz_t}{dt}|_{t=0}=0$.

Second, we prove that $\phi_{d,n}$ is nonsingular by contradiction. Assume that the map $t\to f_t$ is nonsingular at $t=0$, but $t\to f^n_t$ is
singular
at $t=0$. Then for any repelling periodic point $z_0\in \P^1$, $f_0(z_0)$ is a repelling periodic point of $f_0(z)$. Let  $z_t$ be the holomorphic motion
of the periodic point $z_0$. Therefore, $f_t(z_t)$ is the holomorphic motion of the periodic point $f_0(z_0)$. Because $z_t$ and $f_t(z_t)$ are in the repelling cycle of $f_t$ and $t\to f^n_t$ is singular at $t=0$. These
two motions are singular at $t=0$. Then by taking the derivative of $\psi_{1}(z_t,t)=f_t(z_t)$ with respect to $t$,
$$\frac{d\psi_1(z_t, t)}{dt}|_{t=0}=\frac{\partial \psi_1}{\partial z_1}(z_0, 0)\frac{\partial z_t}{\partial t}|_{t=0}+\frac{\partial\psi_1}
{\partial z_2}(z_0,0)=
\frac {df_t(z_t)}{dt}|_{t=0}$$
Since the motions  $z_t$ and $f_t(z_t)$ are singular at $t=0$, i.e. $\frac{dz_t}{dt}|_{t=0}=0$ and $\frac {df_t(z_t)}{dt}|_{t=0}=0$, we can reduce the above equation to
$$\frac{\partial\psi_1}{\partial z_2}(z_0,0)=0.$$
As we know from the previous discussion and the assumption that the map $t\to f_t$ is nonsingular at $t=0$,
$\frac{\partial\psi_1}{\partial z_2}(z,0)=\frac{df_t}{dt}|_{t=0}(z)$ is a nonzero rational function with degree at most $2d$. It has finitely many zeros. However,
the set of repelling periodic points of degree$\geq 2$ rational functions is infinite. This contradicts to the fact that
$\frac{\partial\psi_1}{\partial z_2}(z,0)$ vanishes at any repelling periodic points of $f_0$. \qed

\medskip

{\bf Proof of Theorem \ref{injective}.} As we know,  any regular map from $\P^{2d-1}$ to $\P^{2d^n-1}$ is closed, i.e. the image of any Zariski closed set is Zariski closed. And since the map $\phi_{d,n}: \Rat_d \to \Rat_{d^n}$ is the restriction of a  regular map from $\P^{2d-1}$ to $\P^{2d^n-1}$, the image $\phi_{d,n}(\Rat_d)$ is a subvariety of  $\Rat_{d^n}$. And by Theorem
 5.3 in \cite{H1}, the singularities Sing($\phi_{d,n}(\Rat_d)$) of $\phi_{d,n}(\Rat_d)$ is a proper Zariski closed subset of $\phi_{d,n}(\Rat_d)$. Because the map $\phi_{d,n}$ is
 regular, we have that the preimage $A=\phi_{d,n}^{-1}($Sing($\phi_{d,n}(\Rat_d))$) of Sing($\phi_{d,n}(\Rat_d)$) is a proper
  Zariski closed subset of $\Rat_d$.

Choose a polynomial $p\in \Rat_d$ such that the symmetry group $\Sigma_{J_p}$ of its Julia set is trivial. From the main theorem of
\cite{SS} and Corollary \ref{poly id}, there is no other $f\in \Rat_d$ such that $\phi_{d,n}(f)=\phi_{d,n}(p)$. And since $\phi_{d,n}$ is a proper (preimage of any compact set is compact) regular
map and $\phi_{d,n}$ is an immersion by Proposition \ref{singular}, $p$ is not in $A$.

The set $A\subset \Rat_d$ is $\phi_{d,n}$'s preimage of some Zariski closed subset of $\Rat_{d^n}$, and then $A$ is a proper Zariski closed subset of $\Rat_d$. After throwing away the singularities of $\phi_{d,n}(\Rat_d)$, by Proposition \ref{singular},
$\phi_{d,n}(\Rat_d\backslash A)$ is a connected smooth submanifold of $\Rat_{d,n}$. Moreover, since
 $\phi_{d,n}$ is a proper nonsingular map, the map $\phi_{d,n}$ is a covering map restricted to the following sets:
$$\phi_{d,n}: \Rat_d\backslash A \to \phi_{d,n}(\Rat_d\backslash A).$$
Because $p\notin A$ and there is no other $f\in \Rat_d$ such that $\phi_{d,n}(f)=\phi_{d,n}(p)$. The degree of the covering map should be $1$, i.e.
$\phi_{d,n}: \Rat_d\backslash A \to \phi_{d,n}(\Rat_d\backslash A)$ is injective.

Finally, since $\phi_{d,n}(z^d)=\phi_{d,n}(\zeta z^d)$ for any $\zeta$ with $\zeta^{d^{n-1}+d^{n-2}+\cdots+d^0}=1$, $\phi_{d,n}$ is not injective in $\Rat_d$. Consequently, $A$ is nonempty. \qed

\medskip
Theorem \ref{injective} states that $\phi_{d,n}$ is injective at general points $f\in \Rat_d$. The next theorem indicates that for generic rational functions $f\in \Rat_d$, any rational function $g$ sharing an iterate with $f$ should be some iterate of $f$.

\begin{theorem}\label{aaaa}
For the generic rational functions $f\in \Rat_d$  with degree $d\geq 2$, we have that any rational function $g$,
$$f^n=g^m \Longleftrightarrow g=f^{n/m} \textup{ and } m|n.$$
\end{theorem}
\proof First, for any $d\geq 2$ and $n\geq 2$, $\phi_{d,n}(\Rat_d)$ is a proper subvariety of $\Rat_{d^n}$. Then the set of
rational functions, with degree $d$ and coming from the iteration of lower degree rational functions, is a proper Zariski closed subset of $\Rat_d$. Then for
general $f\in \Rat_d$, it is not an iterate of some lower degree rational function.

Second, for any $d\geq 2$ and $m, n\geq 2$, consider the iteration maps:
$$\Rat_d \to \Rat_{d^n} \to \Rat_{d^{mn}}$$
given by $\phi_{d,n}$ and $\phi_{d^n,m}$. Let $A_{d, mn}\subset \Rat_d$ and $A_{d^n, m}\subset \Rat_{d^n}$ be preimage of the singularities of
$\phi_{d,nm}(\Rat_d)$ and $\phi_{d^n,m}(\Rat_{d^n})$. The set $B=A_{d,nm}\cup \phi^{-1}_{d,n}(A_{d^n, m})$ is a Zariski closed subset of $\Rat_d$. Choose
a polynomial $p\in \Rat_d$ with trivial symmetry group $\Sigma_{J_p}$. From the proof of Theorem \ref{injective}, $p$ is not in $A_{d, nm}$, and $\phi_{d,n}(p)$ is
not in $A_{d^n, m}$. Then  $p\in \Rat_d\backslash B$, i.e. $B$ is a proper Zariski closed subset of $\Rat_d$. And from  the choice of $B$, we know for any $f\in \Rat_d\backslash B$,
$$f^{mn}=g^m \Longleftrightarrow g=f^n.$$

Third, let $d_1, d_2, n_1, n_2\geq 2$ be integers with $d_1^{n_1}=d_2^{n_2}$ and $n_2$ is not divisible by $n_1$. We claim that
$$\phi_{d_1, n_1}(\Rat_{d_1})\nsubseteq \phi_{d_2, n_2}(\Rat_{d_2}).$$
Actually, from the main theorem of \cite{SS} and Corollary \ref{poly id}, there is a polynomial $q\in \Rat_{d_1}$ such that $M_q=\{q, q^2, q^3, \cdots\}$. If there is an
 $h\in \Rat_{d_2}$ such that $h^{n_2}=q^{n_1}$, then $h$ must be in $M_f$. Since $M_q=\{q, q^2, q^3, \cdots\}$, there is some $i$ such that $q^i=h$. Consequently, $q^{n_1}=q^{i*n_2}=h^{n_2}$. So $n_2|n_1$, which contradicts to the assumption that $n_2\nmid n_1$. Consequently, $q^{n_1}\in \phi_{d_1, n_1}(\Rat_{d_1})\nsubseteq \phi_{d_2, n_2}(\Rat_{d_2})$. Then  $\phi^{-1}_{d_1,n_1}(\phi_{d_2, n_2}(\Rat_{d_2}))$
 is a proper Zariski closed subset of $\Rat_{d_1}$. And for any $f\in \Rat_{d_1}\backslash \phi^{-1}_{d_1,n_1}(\phi_{d_2, n_2}(\Rat_{d_2})) $, there is no
 $g\in \Rat_{d_2}$ such that $f^{n_1}=g^{n_2}$.

From the above three statements, we can remove countably many proper Zariski closed subsets of $\Rat_d$.  The left rational functions $f\in \Rat_d$ satisfy the statement of this theorem. \qed

\medskip
With all these preparations, we are ready to prove Theorem \ref{generic rational}.

{\bf Proof of the Theorem \ref{generic rational}.} The first statement of Theorem \ref{generic rational} is just a consequence of Theorem \ref{mainthm}
 and \ref{aaaa}.

Let $f\in \Rat_2$ be rational function with two critical orbits and none of them is preperiodic. Let $g$ be any rational function with $\mu_g=\mu_f$.
Since all  exceptional functions are post-critical finite and $f$ is not post-critical finite, $f$ cannot be exceptional. So by Theorem \ref{main1}, there are integers $m,n,k\geq 1$, such that
$$g^m=(f^{-k}\circ f^k)\circ f^n.$$
By Theorem \ref{genus0} and \ref{3333}, $f^{-k}\circ f^k$ should be  $\sigma_f$ or the identity map. So it indicates
$$g^m=\sigma_f \circ f^n\textup{ or } f^n.$$
By Theorem \ref{aaaa}, there is a generic subset $C\subset \Rat_2$, such that for any $f\in C$,
$f^n=g^m$ implies that $m|n$ and $g=f^{n/m}$, and $\mu_f=\mu_g$ implies that $f^k=\sigma_f\circ g^l$ or $g^l$ for some $k,l\geq 1$.

For any $f(z)=\frac{az^2+bz+c}{dz^2+ez+r}\in \Rat_2$, we can write $\sigma_f$ down explicitly,
$$\sigma_f(z)=\frac{(ar-cd)z-(br-ce)}{(ae-bd)z+(ar-cd)}.$$
There is a free and order two automorphism $\rho$ of $\Rat_2$,
$$\rho: \Rat_2\to \Rat_2,$$
given $\rho(f)=\sigma_f\circ f$.

Since $\rho$ is an automorphism, $C\bigcap \rho^{-1}(C)$ is a generic subset of $\Rat_2$. For any $f\in C\bigcap \rho^{-1}(C)$ and any $g$ with $\mu_f
=\mu_g$, we have $f^n=\sigma_f\circ g^m$ or $g^m$ for some $m,n\geq 1$. If $f^n=g^m$, then $m|n$ and $g=f^{n/m}$. If $f^n=\sigma_f \circ g^m$, then
$g^m=(\sigma_f\circ f)^n$. Since $f \in C\bigcap \rho^{-1}(C)$, it indicates that $\sigma_f\circ f\in C$. So we have $m|n$ and $g=(\sigma_f \circ f)^{n/m}=
\sigma_f \circ f^{n/m}$. \qed

\bigskip

\section{Rational functions with common iterates}\label{common iterates}
In this section, we characterize the condition that two non-exceptional rational functions share an iterate.
\medskip

Theorem \ref{mainthm} says that: generally, having the same measure of maximal entropy is the same as sharing an iterate.  It is easy to see that two rational functions sharing an iterate should have the same set of periodic points. Conversely, for non-exceptional rational functions, having the same set of periodic points also guarantees that they share an iterate. However, this is not true for exceptional functions; see Proposition  \ref{prop5}.

\begin{theorem}[restatement of Theorem \ref{comthm}]
Let $f$ and $g$ be non-exceptional rational functions with degrees $\geq 2$. The following statements are equivalent:
\begin{enumerate}
\item  $f$ and $g$ have some common iterate, i.e. $f^n=g^m$ for some $n, m \in \N^*$.
\item There is some $\varphi$ with degree $\geq 2$, such that $f\circ \varphi= \varphi \circ f$ and $g\circ \varphi=\varphi \circ g$.
\item The maximal entropy measures $\mu_f=\mu_g$, and $J\cap \textup{Per}(f)\cap \textup{Per}(g)\neq \emptyset$.
\item $\Prep(f)=\Prep(g)$ and $J\cap \textup{Per}(f)\cap \textup{Per}(g)\neq \emptyset$.
\item $ \textup{Per}(f)=\textup{Per}(g)$.
\end{enumerate}
\end{theorem}
 {\bf Proof Theorem \ref{comthm}.} For (1), let $\varphi=f^n=g^m$. Since $f$ and $g$ are both commutable with $\varphi$, (1) $\Rightarrow$ (2). By Ritt's theorem, non-exceptional commutable rational functions share an iterate; see \cite{R1} and also \cite{E1}. Since $f$ and $g$ are non-exceptional, $f\circ \varphi= \varphi \circ f$ and $g\circ \varphi=\varphi \circ g$ implies that $\phi$ is non-exceptional, and $f, g, \phi$ share an iterate. So (2) $\Rightarrow$ (1).

By Yuan and Zhang's Theorem 1.6 in \cite{YZ}, PrePer($f$)$=$PrePer($g$) implies that $\mu_f=\mu_g$. Then we have (4)$\Rightarrow$(3). Two rational functions sharing an iterate must have the same set of periodic (preperiodic) points  and also the same measure of maximal entropy. It has (1)$\Rightarrow$(4) and (1)$\Rightarrow$(5).

From Theorem 3.5 in \cite{YZ} and Theorem 1.2 in \cite{BD}, it shows that two rational functions  sharing infinitely many preperiodic points guarantees they have common set of periodic points. Assume that Per($f$)$=$Per($g$). Since Per($f$) is not a finite set,  PrePer($f$)$=$PrePer($g$). It shows that (5)$\Rightarrow$(4).

It remains to show that (3)$\Rightarrow$ (1). The proof is inspired by \cite{L1}, which follows from the description of local dynamics in \cite{L4}. Suppose $\mu_f=\mu_g$, and $J\cap \textup{Per}(f)\cap \textup{Per}(g)\neq \emptyset$. By passing to some iterates and changing  of coordinates, we can assume that $0$ is a fixed point of $f$ and $g$, and $0$ is in  their Julia set. Since $0$ is in the Julia set and it is fixed by both $f$ and $g$, then  $R(z)=f^{-1}\circ g^{-1} \circ f \circ g(z)$ is locally well defined near $0$.
We claim that $R$ is the identity map. Otherwise, since $R$ has multiplier equaling $1$ at its fixed point $0$. It determines attracting and repelling flowers near $0$. Suppose that there is some point $x$ near $0$ in the Julia set and also in some attracting petal of the flowers
determined by $R$; see Section 10 of \cite{M2}. Then there is some fundamental domain of $R$ for this petal, which contains some neighborhood of this point $x$. As $\mu_f$ is supported in the Julia set, the fundamental domain won't have zero measure.  Since $R$ acts on this petal like a
transformation (in appropriate coordinate) and $R$ preserves the measure $\mu_f$, the $\mu_f$-measure of this petal  cannot be finite. However, we know that the total mass of $\mu_f$ on $\P^1$ is $1$, which is a contradiction. So there is no point in the Julia set which is in the attracting petals of the flowers of
$R$. Replace $R$ by $R^{-1}$, similarly we  see that there is no point in the Julia set which is in the attracting petals of the  flowers of $R^{-1}$. As the union of the attracting petals of $R$ and $R^{-1}$ contains a small disc punctured at $0$. The point  $0$ should be an isolated point
 in the Julia set of $f$ and $g$. This is impossible, since a Julia set cannot have isolated points. Therefore, $R$ should be the identity map. Which means $f$ and $g$ are commutable. So the third statement implies the first  statement.

 \qed

For exceptional case, we would not have this nice result.

\begin{prop}\label{prop5}
Let $f(z)=z^{d_f}$ and $g(z)=z^{d_g}$, with $d_f, d_g\geq 2$. Then
$$\Per (f)=\Per (g)\Leftrightarrow \textup{$\forall$ prime $p$, $p|d_f$ iff $p|d_g$}.$$
\end{prop}
\proof Assume $p\geq 2$ is a prime number such that $p|d_f$ and $p\nmid d_g$. There are integers $n\geq 1$ and $m$, with
$$d_g^n=mp+1$$
Let $z_o=e^{2\pi i/p}$, we have $f(z_o)=1$ and $g^n(z_o)=e^{2\pi i (mp+1)/p}=z_o$. So $z_o$ is preperiodic and not periodic for $f$ but it is  periodic for $g$, i.e. $\Per (f)\neq\Per (g)$.

Conversely, assume that for any prime number $p$, $p|d_f$ iff $p|d_g$. Let $z_o=e^{2a\pi i/b}$ be a periodic point of $f$ with period $n$, where $a$ and $b$ are coprime integers. Then we have $d_f^na/b=a/b+m$ for some integer $m$.
$$d_f^na=a+mb\Rightarrow (d_f^n-1)a=mb$$
So $b|(d_f^n-1)$, which means that $b$ and $d_f$ are coprime integers. Then by previous assumption, $b$ and $d_g$ are coprime integers. There is some integer $k\geq 1$ such that $d_g^k\equiv1$(mod $b$), i.e. $d_g^ka/b=a/b+t$ for some integer $t$. In all we have $g^k(z_o)=z_o$. So $z_o$ is a periodic point of $g$. Consequently, Per($f$)$\subset$Per($g$). Similarly, we have Per($g$)$\subset$Per($f$). In all, Per($g$)$=$Per($f$).
\qed

\medskip
Let $f(z)=z^{2*3}$ and $g(z)=z^{4*3}$. From the above proposition, they have the same set of periodic points. However, they do not share an iterate, which can be seen from the degrees of them.

There is one more thing  I would like to mention here. As we know, for any two rational functions, if the intersection of their sets of preperiodic points has
infinitely many points, then they have identical set of preperiodic points. For any two non-exceptional rational functions $f$ and $g$, $|$Per($f$)$\bigcap$Per($g$)$|$$=\infty$ guarantees that Per($f$)$=$Per($g$).  However, for exceptional polynomials $f(z)=z^3$ and $g(z)=z^5$, we have $|\Per(f)\cap\Per(g)|=\infty$ and $\Per(f)\neq \Per(g)$ by Proposition \ref{prop5}, since it is not hard to see that $e^{2\pi i/2^k}$ is a periodic point for both $f$ and $g$ with any $k\geq 1$.

%%%%%%
%%%%%%
%%%%%%
%%%%%%

\bigskip\bigskip
\def\cprime{$'$}

%\bibliographystyle{../tex/bib/math}
%\bibliography{../tex/bib/math}
%\bibliographystyle{math}
%\bibliography{refs}

 \end{document}